\documentclass[12pt]{amsart}
\date{\today}

\usepackage[cp1251]{inputenc}         

\usepackage[ukrainian]{babel} 
\usepackage{amsmath,amsthm,amsfonts,amssymb,mathrsfs}

\usepackage{color}

\usepackage{hyperref}

\newtheorem{theorem}{Теорема}

\newtheorem{proposition}[theorem]{Твердження}
\newtheorem{corollary}[theorem]{Наслiдок}
\newtheorem{lemma}[theorem]{Лема}
\theoremstyle{definition}
\newtheorem{example}[theorem]{Приклад}
\newtheorem{remark}[theorem]{Зауваження}

\begin{document}

\title[Про локально компактнi групи з нулем]{Про локально компактнi групи з нулем}

\author[Катерина Максимик]{Катерина Максимик}
\address{Механіко-математичний факультет, Львівський національний університет ім. Івана Франка, Університецька 1, Львів, 79000, Україна}
\email{kate.maksymyk15@gmail.com}



\begin{abstract}
В роботі досліджуються  алгебраїчні умови на групу $G$, при виконанні яких локально компактна трансляційно неперервна топологія на дискретній групі $G$ з приєднаним нулем є або компактною, або дискретною. Введено  електорально гнучкі та електорально стійкі групи та вивчаються їх властивості. Зокрема, доведено, що кожна група, яка  містить нескінченну циклічну підгрупу нескінченного індексу та кожна незліченна комутативна група є електорально гнучкими, а також, що кожна зліченна локально скінченна група є електорально стійкою. Основним результатом роботи є таке твердження: якщо $G$~--- дискретна електорально гнучка нескінченна група, то  кожна гаусдорфова трансляційно неперервна локально компактна топологія на $G^0$ є або дискретною, або компактною. На довільній нескінченній віртуально циклічній групі (а отже, на електорально стійкій групі)   з приєднаним нулем $G^0$ побудовано недискретну некомпактну локально компактну трансляційно неперервну топологію, яка індукує на $G$ дискретну топологію.

\smallskip 

\emph{Kateryna Maksymyk}, \textbf{On locally compact groups with zero}.

We study algebraic properties on a group $G$ such that if the discrete group $G$ has these properties then every locally  compact shift continuous topology on $G$ with adjoined zero is either compact, or discrete. We introduce  electorally flexible and electorally stable groups and establish their properties. In particular, we prove that every group with an infinite cyclic subgroup of an infinite index and every uncountable commutative group are electorally flexible, and show that every countable locally finite group is electorally stable. The main result of the paper is the following: if $G$ is a discrete electorally flexible group then every Hausdorff locally compact shift-continuous topology on $G$ with adjoined zero is either compact, or discrete. Also, we construct a non-discrete non-compact Hausdorff locally compact shift-continuous topology on any discrete virtually cyclic group (and hence on a electorally stable group)  $G$ with adjoined zero.
\end{abstract}

\maketitle



У даній праці ми користуватимемося термінологією з \cite{Kurosh-1967, Clifford-Preston-1961-1967, Ruppert-1984}. Усі простори вважаються гаусдорфовими, якщо не зазначено інше.

Надалі для групи $G$ через $G^0$ позначається група $G$ з приєднаним нулем \cite{Clifford-Preston-1961-1967}. У комутативній групі $G$ добуток двох елементів $a$ i $b$ позначатимемо через $a+b$, а обернений елемент до $a\in G$ позначатимемо через $-a$.

Напівгрупа $S$ називається \emph{інверсною}, якщо для довільного елемента $x\in S$ існує єдиний елемент $y\in S$ такий, що $xyx=x$ i $yxy=y$. У цьому випадку кажуть, що $y$ є \emph{інверсним} елементом до $x$ в $S$, а відображення з $S$ у $S$, що ставить кожному елементу його інверсний, називається \emph{інверсією} \cite{Clifford-Preston-1961-1967}. Очевидно, що група з приєднаним нулем є інверсною напівгрупою.

Напівгрупа $S$ із заданою на ній топологією $\tau$ називається (\emph{напів})\emph{топологічною}, якщо напівгрупова операція в $(S,\tau)$ є (нарізно) неперервною, і в цьому випадку кажуть, що $\tau$ є \emph{напівгруповою} (\emph{трансляційно неперервною}) топологією на $S$ \cite{Carruth-Hildebrant-Koch-1983-1986}. Інверсна топологічна напівгрупа $(S,\tau)$ з неперервною інверсією називається \emph{топологічною інверсною напівгрупою}, а топологія $\tau$ в цьому випадку називається \emph{інверсною напівгруповою} топологією на $S$.

Відомо, якщо $S$~--- напівтопологічна напівгрупа з нулем $0$ така, що $S\setminus \{0\}$~--- компактний простір, чи $G^0$ --- група з приєднаним нулем, що є топологічною інверсною напівгрупою, то нуль $0$ є ізольованою точкою. В загальному випадку навіть для локально компактних груп з приєднаним нулем, які є топологічними напівгрупами, це не так \cite{Carruth-Hildebrant-Koch-1983-1986, Koch-Wallace-1964}. Однак нуль в компактній топологічній $0$-простій напівгрупі є ізольованою точкою \cite{Paalman-de-Miranda-1964}. Проте ці результати на можна поширити на локально компактні цілком $0$-прості топологічні напівгрупи \cite{Owen-1973} і на компактні цілком $0$-прості напівтопологічні зліченні напівгрупи \cite{Gutik-Pavlyk-2005}. Карл Гофманн у праці \cite{Hofmann-1963} описав структуру локально компактної топологічної групи з приєднаним нулем у випадку локально компактних топологічних напівгруп. Приєднання нуля до близьких до компактних (напів)топологічних груп та топологія в нулі для деяких класів локально компактних напівгруп вивчалось в працях \cite{Demenchuk-1986, Bardyla-2018, Bardyla-Gutik-2016, Gutik-2015, Gutik-2018, Gutik-Ravsky-2015, Mesyan-Mitchell-Morayne-Peresse-2016, Owen-1973, Roset-1977, Shneperman-1981, Shneperman-1985}.

У даній праці досліжуються алгебраїчні умови на групи, при виконанні яких локально компактна трансляційно неперервна топологія на дискретній групі з приєднаним нулем є або дискретною, або компактною.

Дослідження цієї праці мотивовані проблемою 7 з \cite{Berglund-1980}: ``\emph{яка хороша структура замикання групи в напівтопологічній напівгрупі?}'' та наступним простим твердженням.

\begin{proposition}\label{proposition-1}
Якщо $T_1$-напівтопологічна напівгрупа $S$ містить власну щільну дискретну підгрупу $G$, то $I=S\setminus G$~--- двобічний ідеал в $S$.
\end{proposition}

\begin{proof}
За лемою~3 з \cite{Gutik-Savchuk-2017}, $G$ --- відкритий підпростір в $S$.

Зафіксуємо довільний елемент $y\in I$. Якщо $xy=z\notin I$ для деякого $x\in G$, то існує відкритий окіл $U(y)$ точки $y$ у просторі $S$ такий, що $\{x\}\cdot U(y)=\{z\}\subset G$. Окіл $U(y)$ містить нескінченну кількість елементів групи  $G$, а це суперечить тому, що в групі зсуви є бієктивними відображеннями. З отриманого протиріччя випливає, що $xy\in I$ для всіх $x\in G$ й $y\in I$. Доведення того, що  $yx\in I$ для всіх $x\in G$ і $y\in I$ а аналогічним.

Припустимо протилежне: нехай $xy=w\notin I$ для деяких $x,y\in I$. Тоді $w\in G$ і з нарізної неперервності напівгрупової операції в $S$ випливає, що існують відкриті околи  $U(x)$ і $U(y)$ точок $x$ і $y$ в $S$, відповідно, такі, що $\{x\}\cdot U(y)=\{w\}$ і $U(x)\cdot \{y\}=\{w\}$. Оскільки обидва околи $U(x)$ і $U(y)$ містять нескінченну кількість елементів групи $G$, то кожна з рівностей $\{x\}\cdot U(y)=\{w\}$ і $U(x)\cdot \{y\}=\{w\}$ суперечить тому, що в групі зсуви є бієктивними відображеннями. З отриманого протиріччя випливає, що $xy\in I$.
\end{proof}

Будемо говорити, що нескінченна група $G$ є:\footnote{Електорально гнучкі та електорально стійкі групи введені проф. Т. О. Банахом, і були анонсовані на семінарі ``Топологія та її застосування'' у Львівському університеті в 2019 році.}
\begin{itemize}
  \item \emph{електорально гнучкою}, якщо для довільного розбиття $G=A\sqcup B$ групи $G$ на дві нескінченні множини, існують нескінченна множина $I\subseteq A$ та елемент $x\in G$ такі, що $I\cdot x\subseteq B$;
  \item \emph{електорально стійкою}, якщо $G$ не є електорально гнучкою.
\end{itemize}

Підмножина $A$ групи $G$ називається \emph{трансляційно майже стійкою}, якщо для довільного $x\in G$ симетрична різниця  $A\Delta (A\cdot x)$ є скінченною. Очевидно, що кожна скінченна підмножина в довільній групі, а також коскінченні підмножини в нескінченних групах є трансляційно майже стійкими. Також, в адитивній групі цілих чисел множина натуральних чисел є трансляційно майже стійкою. Виникає природне питання: \emph{у яких нескінченних групах їх трансляційно майже стійкі підмножини вичерпуються скінченними та коскінченними підмножинами?} Це питання також мотивоване наступною дихотомією локально компактних напівтопологічних груп з приєднаним нулем:

\begin{lemma}\label{lemma-2}
Нехай $G$~--- дискретна група така, що кожна нескінченна трансляційно майже стійка підмножина в $G$ є коскінченною. Тоді кожна гаусдорфова трансляційно неперервна локально компактна топологія на $G^0$ є або дискретною, або компактною.
\end{lemma}

\begin{proof}
Очевидно, що дискретна топологія на $G^0$ є трансляційно неперервною та локально компактною.

Нехай $\tau$~--- недискретна гаусдорфова трансляційно неперервна локально компактна топологія на $G^0$ i $U_0$~--- нескінченний відкритий компактний окіл нуля $0$ в просторі $(G^0,\tau)$. Оскільки кожен лівий чи правий зсув у $(G^0,\tau)$ на елемент групи $G$ є гомеоморфізмом, то симетрична різниця $U_0\Delta (U_0\cdot x)$ є скінченною для довільного $x\in G$. З припущення теореми випливає, що $U_0\setminus\{0\}$~--- коскінченна підмножина в $G$.
\end{proof}

\begin{proposition}\label{proposition-3}
Нескінченна група $G$ є електорально гнучкою тоді і лише тоді, коли кожна нескінченна трансляційно майже стійка підмножина в $G$ є коскінченною.
\end{proposition}

\begin{proof}
Припустимо, що існує електорально гнучка нескінченна група $G$, що містить нескінченну трансляційно майже стійку підмножину $A$ в $G$ з нескінченним доповненням $B=G\setminus A$. Тоді існують нескінченна множина $I\subseteq A$ й елемент $x\in G$ такі, що $I\cdot x\subseteq B$, а це суперечить трансляційній майже стійкості підмножини $A$ в $G$.

Припустимо, що в групі $G$ кожна нескінченна трансляційно майже стійка підмножина  є коскінченною, але група $G$ не є електорально гнучкою. Тоді існує розбиття $G=A\sqcup B$ групи $G$ на дві нескінченні множини таке, що для довільних нескінченної множини $I\subseteq A$ та елемента $x\in G$ такі, що множина $I\cdot x\cap B$ скінченна. Отже,  $A\cdot x\cap B=A\cdot x\cap(G\setminus A)$ --- скінченна множина, а це суперечить нашому припущенню.
\end{proof}

З леми~\ref{lemma-2} і твердження~\ref{proposition-3} випливає

\begin{theorem}\label{theorem-4}
Нехай $G$~--- дискретна електорально гнучка нескінченна група. Тоді кожна гаусдорфова трансляційно неперервна локально компактна топологія на $G^0$ є або дискретною, або компактною.
\end{theorem}

\begin{corollary}\label{corollary-5}
Нехай $G$~--- дискретна електорально гнучка нескінченна група. Тоді кожна гаусдорфова напівгрупова локально компактна топологія на $G^0$ є  дискретною.
\end{corollary}

\begin{proposition}\label{proposition-6}
Нехай $G$~---  електорально гнучка нескінченна зліченна група. Тоді кожна гаусдорфова трансляційно неперервна локально компактна топологія на $G^0$ є або дискретною, або компактною.
\end{proposition}

\begin{proof}
Нехай $\tau$~--- гаусдорфова трансляційно неперервна локально компактна топологія на $G^0$. Оскільки нуль напівгрупи $G^0$ є замкненою підмножиною в $(G^0,\tau)$, то за наслідком~3.3.10 з \cite{Engelking-1989}, $G$~--- локально компактний простір. Оскільки група $G$ --- зліченна, то за теоремою Бера про категорії (див. теорема~3.9.3 з \cite{Engelking-1989}) простір $G$ містить ізольовану в $G$, а отже і в $(G^0,\tau)$, точку. З того, що всі зсуви в $(G^0,\tau)$ на елементи групи $G$ є гомеоморфізмами випливає, що всі точки в $G$ є ізольованими. Далі скористаємося теоремою~\ref{theorem-4}.
\end{proof}

\begin{corollary}\label{corollary-7}
Нехай $G$~---  електорально гнучка нескінченна зліченна група. Тоді кожна гаусдорфова напівгрупова локально компактна топологія на $G^0$ є   дискретною.
\end{corollary}

Наступні два твердження дають достатні умови, при виконанні яких група є електорально гнучкою.

Нагадаємо \cite{Kurosh-1967}, що \emph{індексом підгрупи} ${\displaystyle H}$ у групі ${\displaystyle G}$ називається потужність множини класів суміжності в кожному (правому або лівому) із розкладів групи ${\displaystyle G}$ за цією підгрупою ${\displaystyle H}$.

\begin{proposition}\label{proposition-8}
Якщо група $G$ містить нескінченну циклічну підгрупу $Z\subset G$ нескінченного індексу, то $G$ є електорально гнучкою.
\end{proposition}

\begin{proof}
Нехай $z$~--- породжуючий елемент групи $Z$. Розглянемо розбиття $G=A\sqcup B$ групи $G$ на дві нескінченні множини. Розглянемо множини
\begin{equation*}
J_+=\left\{a\in A\colon a\cdot z\in B\right\} \qquad \hbox{i} \qquad J_-=\left\{a\in A\colon a\cdot z^{-1}\in B\right\}.
\end{equation*}
Якщо одна з цих множин є нескінченною, то доведення завершено.  Тому, ми припустимо що множина $J=J_+\cup J_-$ є скінченною.

Якщо множини $A\setminus (J\cdot Z)$ i $B\setminus (J\cdot Z)$ є непорожніми, то ми можемо зафіксувати точки $a\in A\setminus (J\cdot Z)$ i $b\in  B\setminus (J\cdot Z)$, і зауважимо, що для нескінченної множини $I=a\cdot Z\subseteq A$ і точки $x=ba^{-1}$ маємо
\begin{equation*}
x\cdot I=b\cdot a^{-1}\cdot a \cdot Z = b\cdot Z\subseteq B.
\end{equation*}

Отже, ми можемо припускати, що одна з множин $A\setminus (J\cdot Z)$ або $B\setminus (J\cdot Z)$ є порожньою. Якщо множина $A\setminus (J\cdot Z)$ є порожньою, то $A\subseteq J\cdot Z$ і за принципом Діріхле (див. \cite[підрозділ~3.1]{Brualdi-2009}) для деякого елемента $a\in A$ множина $I=A\cap (a\cdot Z)$ є нескінченною. Тоді для довільного елемента $b\in G \setminus (J\cdot Z)\subseteq B$ маємо, що $b\cdot a^{-1}\cdot I\subseteq b\cdot Z\subseteq B$.

Якщо множина $B\setminus (J\cdot Z)$ є порожньою, то $B\subseteq J\cdot Z$ і для деякого елемента $b\in B$ множина $B\cap (b\cdot Z)$ є нескінченною. Тоді для довільного елемента $a\in G\setminus(J\cdot Z)\subseteq A$ маємо, що множина \begin{equation*}
I=a\cdot b^{-1}\cdot(B\cap(b\cdot Z))\subseteq a\cdot A\subseteq A
\end{equation*}
є нескінченною та $b\cdot a^{-1}\cdot I\subseteq B$.
\end{proof}

Нагадаємо \cite{Robinson-1996}, що група $G$ називається \emph{локально скінченно}ю, якщо кожна її скінченна підмножина міститься в скінченній підгрупі в $G$.

\begin{proposition}\label{proposition-9}
Кожна незліченна комутативна група $G$ є електорально гнучкою.
\end{proposition}

\begin{proof}
Якщо група $G$ не є локально скінченною, то $G$ містить нескінченну циклічну підгрупу незліченного індексу, і тоді за твердженням~\ref{proposition-8} група $G$ є електорально гнучкою. Отже, ми припускатимемо, шо група $G$ є локально скінченною. Нехай $G=A\sqcup B$~--- розбиття групи $G$ на дві нескінченні множини. Якщо множина $A$ є зліченною, то виберемо довільний елемент $x\in G\setminus A-A$ i стверджуємо, що $(x+A)\cap A=\varnothing$, а отже $x+A\subseteq B$. Якщо множина $B$ є зліченною, то можемо вибрати елемент $x\in G$ такий, що $I= (x+B)\subseteq A$, а отже $-x+I\subseteq B$.

Далі ми припускатимемо, що обидві множини $A$ i $B$ є незліченними. Зафіксуємо довільну нескінченну зліченну підгрупу $Z\subset G$. Якщо для деякого елемента $z\in Z$ множина $(A+z)\cap B$ є нескінченною, то доведення завершене. Отже, ми припускатимемо, що для кожного елемента $z\in Z$ множина $F_z=(A+z)\cap B$ є скінченною. Тоді множина
\begin{equation*}
  S=(A+Z)\cap (B+Z)=\bigcup_{z,z'\in Z}(A+z)\cap(B+z')
\end{equation*}
є непорожньою і щонайбільше зліченною. Позаяк множини $A$ i $B$ є незліченними, то існують елементи $a\in A\setminus S$ i $b\in B\setminus S$. Тоді для нескінченної множини $I=a+Z\subseteq A$ та елемента $x=b-a$ отримуємо, що $x+I=b+Z\subseteq B$, звідки випливає, що група $G$ є електорально гнучкою.
\end{proof}

\begin{proposition}\label{proposition-10}
Кожна зліченна локально скінченна група $G$ є електорально стійкою.
\end{proposition}

\begin{proof}
Зобразимо групу $G$ як об'єднання $\bigcup_{n\in\omega}G_n$ строго зростаючої послідовності скінченних груп. Для довільного $n\in\omega$ зафіксуємо підмножину $A_n \subseteq G_{n+1} \setminus G_n$, яка має одноточковий перетин $A_n\cap (x\cdot G_n)$ з кожним суміжним класом $x\cdot G_n$ та елементом  $x\in G_{n+1}\setminus G_n$.

Приймемо
\begin{equation*}
A= G_0\cup \bigcup_{n\in\omega}(A_{2n}\cdot G_{2n}) \qquad \hbox{i} \qquad B =\bigcup_{n\in\omega}(A_{2n+1}\cdot G_{2n+1}).
\end{equation*}
Тоді для довільного елемента $x\in G$ множина $(A\cdot x)\cap B$ є скінченною.
\end{proof}

Нагадаємо \cite{Epstein-1962, Wall-1967}, що група $G$ називається \emph{віртуально циклічною}, якщо $G$ містить циклічну підгрупу скінченного індексу.

З твердження~\ref{proposition-10} випливають такі два наслідки.

\begin{corollary}\label{corollary-11}
Нетривіальна комутативна група без кручень $G$ є електорально стійкою тоді і лише тоді, коли $G$ ізоморфна адитивній групі цілих чисел $\mathbb{Z}$.
\end{corollary}

\begin{corollary}\label{corollary-12}
Комутативна група $G$ є електорально стійкою тоді і лише тоді, коли $G$ є або зліченною та локально скінченною, або є віртуально циклічною.
\end{corollary}

Нагадаємо \cite{Epstein-1962, Wall-1967}, що група $G$ має \emph{більше ніж один кінець}, якщо існує нескінченна підмножина $S\subset G$ з нескінченним у $G$ доповненням така, що симетрична різниця $S\Delta (S\cdot x)$ є скінченною для довільного $x\in G$. Будемо також говорити, що група $G$ має \emph{один кінець}, якщо існує така єдина нескінченна підмножина $S\subset G$ із скінченним доповненням, яка задовольняє вищезгадані умови. Очевидно, що нескінченні електорально гнучкі групи --- це в точності групи з одним кінцем.

\medskip

У класичній комбінаторній теорії груп добре відома наступна теорема, яка описує нескінченні групи з двома кінцями:

\begin{theorem}[\cite{Epstein-1962, Wall-1967}]\label{theorem-13}
Для групи $G$ наступні умови є еквівалентними:
\begin{itemize}
  \item[$(i)$] $G$ --- група з двома кінцями;
  \item[$(ii)$] $G$ --- нескінченна віртуально циклічна група.
\end{itemize}
\end{theorem}

\begin{remark}\label{remark-14}
\begin{enumerate}
  \item[1.] З теореми~\ref{theorem-13} випливає, що адитивна група цілих чисел $\mathbb{Z}$, а також її прямий добуток з довільною скінченною групою, має два кінці. Також, на адитивній групі цілих чисел з приєднаним нулем $\mathbb{Z}^\mathbf{0}$ існує рівно чотири гаусдорфові локально компактні трансляційно неперервні топології (див. твердження~4.5 з \cite{Gutik-2018}), причому три з них є напівгруповими.
  \item[2.] З твердження~\ref{proposition-9} випливає, що наступні групи є електорально гнучкими:
  \begin{enumerate}
    \item $n$-а пряма степінь адитивної групи цілих чисел $\mathbb{Z}^n$ для $n\geqslant 2$;
    \item $n$-а пряма степінь адитивної групи раціональних чисел $\mathbb{Q}^n$ для $n\geqslant 1$;
    \item $n$-а пряма степінь адитивної групи дійсних чисел $\mathbb{R}^n$ для $n\geqslant 1$;
    \item $n$-а пряма степінь мультиплікативної  групи дійсних чисел $(\mathbb{C}^*)^n$ для $n\geqslant 1$ та її підгрупа $n$-вимірний тор $\mathbb{T}^n$, як $n$-а пряма степінь одиничного кола  $\mathbb{S}^1=\left\{z\in\mathbb{C}\colon |z|=1\right\}$.
  \end{enumerate}
   \item[3.] З твердження~\ref{proposition-8} випливає, що вільна (абелева) група $F_X$ над множиною $X$ потужності $\geqslant 2$ є електорально гнучкою.
\end{enumerate}
\end{remark}

З наступного прикладу випливає, що на нескінченній віртуально циклічній групі з приєднаним нулем $G^0$ існують недискретні некомпактні локально компактні трансляційно неперервні топології, які індукують на $G$ дискретну топологію.

\begin{example}\label{example-15}
Нехай $G$ --- нескінченна віртуально циклічна група та $K_1$ i $K_2$~--- кінці в групі $G$. Для $i=1,2$ означимо на $G^0$ топологію $\tau_i$ наступним чином:
\begin{enumerate}
  \item усі елементи групи $G$ є ізольованими точками в $(G^0,\tau_i)$;
  \item сім'я $\mathscr{B}_i(0)=\left\{g_1 K_i g_2\cup\{0\}\colon g_1,g_2\in G\right\}$ є базою топології $\tau_i$ в нулі напівгрупи $G^0$.
\end{enumerate}
Оскільки $K_i$~--- кінець в $G$, то $\tau_i$ --- трансляційно неперервна гаусдорфова локально компактна топологія на напівгрупі $G^0$, яка не є ні дискретною, ані компактною.
\end{example}

Зауважимо, що топології $\tau_1$ і $\tau_2$ у випадку адитивної групи цілих чисел $\mathbb{Z}$ є напівгруповими \cite{Gutik-2018}.
\section*{Подяка}

Автор висловлює подяку проф. Т. О. Банаху та науковому керівнику О. В. Гутіку 
за корисні поради, коментарі та  зауваження.


\end{document}